\documentclass{siamart190516}

\usepackage{amssymb,amsfonts,amsmath}

\usepackage[mathscr]{eucal}

\usepackage{stmaryrd}

\usepackage{amsbsy}

\usepackage{bm}

\usepackage{braket}

\usepackage{xspace}

\usepackage{booktabs}

\usepackage{xparse}

\usepackage{cleveref}

\usepackage{algpseudocode}

\Crefname{ALC@unique}{Line}{Lines}

\usepackage{lmodern}
\usepackage{listings}
\usepackage[most]{tcolorbox}

\DeclareTCBListing{macrobox}{s G{} }{IfBooleanTF={#1}{}{listing side text},title=#2,
  listing options={style=tcblatex,commentstyle=\color{red!70!black}}
}

\usepackage{fixme}
\fxsetup{
  status=draft,
  nomargin,
  inline,
  theme=color,
}
\fxsetup{marginface=\linespread{1}\footnotesize}
\FXRegisterAuthor{tk}{tke}{TK}
\FXRegisterAuthor{nj}{nje}{NJ}
\FXRegisterAuthor{ep}{epe}{EP}

\usepackage{subcaption}

\NewDocumentCommand{\qtext} {m} {\quad\text{#1}\quad}

\NewDocumentCommand{\Real} {} {\mathbb{R}}
\NewDocumentCommand{\E} {} {\mathbb{E}}

\NewDocumentCommand{\VcRoot}{m m G{\BooleanFalse} G{\BooleanFalse}}{%
    \bm{#1\mathbf{\MakeLowercase{#2}}}%
    \IfBooleanT{#3}{^{\intercal}}%
    \IfBooleanT{#4}{^{\phantom{\intercal}}}%
  }
\NewDocumentCommand{\Vc}{ O{} m t' t"} {\VcRoot{#1}{#2}{#3}{#4}}

\NewDocumentCommand{\MxRoot}{m m G{\BooleanFalse} G{\BooleanFalse}}{%
    \bm{#1\mathbf{\MakeUppercase{#2}}}%
    \IfBooleanT{#3}{^{\intercal}}%
    \IfBooleanT{#4}{^{\phantom{\intercal}}}%
  }
\NewDocumentCommand{\Mx}{O{} m t' t"}{\MxRoot{#1}{#2}{#3}{#4}}

\NewDocumentCommand{\Tn}{O{} m}{\boldsymbol{#1\mathscr{\MakeUppercase{#2}}}}

\NewDocumentCommand{\Tm}{O{} m m t' t"}{\MxRoot{#1}{#2}{#4}{#5}_{(#3)}}

\NewDocumentCommand{\X} {} {\Tn{X}}

\NewDocumentCommand{\Xk} { O{} G{n} t' t"} {\MxRoot{#1}{X}{#3}{#4}_{#2}}

\DeclareDocumentCommand{\xi}{ s } 
{
  \IfBooleanTF{#1}
  {x_{i_1 \dots i_d}}
  {x_{i}}
}

\def\lastiter{\bar}

\NewDocumentCommand{\Xt} {s} {\IfBooleanTF{#1}{\Mx{X}_t}{\Tn{X}_t}}
\NewDocumentCommand{\Xh} {s} {\IfBooleanTF{#1}{\Mx{X}_h}{\Tn{X}_h}}
\NewDocumentCommand{\xit}{}{x_{it}}
\NewDocumentCommand{\xih}{}{x_{ih}}
\NewDocumentCommand{\yit}{s}{\IfBooleanT{#1}{\tilde}y_{it}}

\NewDocumentCommand{\yiht}{s}{\IfBooleanT{#1}{\tilde}y_{ih}}
\NewDocumentCommand{\cit}{s}{\IfBooleanT{#1}{\tilde}c_{it}}

\NewDocumentCommand{\Wt} {s} {\IfBooleanTF{#1}{\Mx{W}_t}{\Tn{W}_t}}
\NewDocumentCommand{\Wh} {s} {\IfBooleanTF{#1}{\Mx{W}_h}{\Tn{W}_h}}

\NewDocumentCommand{\M} {} {\Tn{M}}
\NewDocumentCommand{\Mt} {s} {\IfBooleanTF{#1}{\Mx{M}_t}{\Tn{M}_t}}
\NewDocumentCommand{\Mtold} {s} {\IfBooleanTF{#1}{\Mx[\lastiter]{M}_t}{\Tn[\lastiter]{M}_t}}
\NewDocumentCommand{\Mh} {s} {\IfBooleanTF{#1}{\Mx{M}_h}{\Tn{M}_h}}
\NewDocumentCommand{\Mhold} {} {\Tn[\lastiter]{M}_{h}}

\DeclareDocumentCommand{\mit}{}{m_{it}}

\NewDocumentCommand{\miht}{}{m_{ih}}
\NewDocumentCommand{\mihold}{}{\lastiter{m}_{ih}}

\NewDocumentCommand{\A}{t'}{\mathbf{A}\IfBooleanT{#1}{^{\mkern -4mu \intercal}}}
\NewDocumentCommand{\B}{}{\Mx{B}}

\NewDocumentCommand{\Ak} { G{k} } {%
  \Mx{A}^{\mkern -4mu (#1)}
}
\NewDocumentCommand{\Akt} { G{k}G{t} } {%
  \Mx{A}^{\mkern -4mu (#1)}_{#2}
}

\NewDocumentCommand{\akj}{O{k} G{j}}{%
  \Vc{a}^{\mkern -3mu {(#1)}}_{#2}
}

\NewDocumentCommand{\Akold} { G{k} } {%
  \Mx[\lastiter]{A}^{\mkern -4mu (#1)}
}

\NewDocumentCommand{\lj}{G{j} t"}{s_{#1}\IfBooleanT{#2}{^{\phantom{(1)}}\!\!}}
\NewDocumentCommand{\lvec}{}{\Vc{s}}

\NewDocumentCommand{\st}{G{t} t'}{\Vc{s}_{#1} \IfBooleanT{#2}{^{\intercal}}}
\NewDocumentCommand{\sh}{}{\Vc{s}_h}

\NewDocumentCommand{\HW}{}{\mathcal{H}}
\NewDocumentCommand{\HWt}{ G{t} }{\mathcal{H}_{#1}}

\NewDocumentCommand{\Y}{s}{\Tn[\IfBooleanT{#1}{\tilde}]{Y}}
\NewDocumentCommand{\Yt}{s}{\Tn[\IfBooleanT{#1}{\tilde}]{Y}_t}
\NewDocumentCommand{\Yh}{s}{\Tn[\IfBooleanT{#1}{\tilde}]{Y}_h}
\NewDocumentCommand{\Yht}{s}{\Tn[\IfBooleanT{#1}{\tilde}]{Y}_{h}}
\NewDocumentCommand{\Yk}{s}{\Mx[\IfBooleanT{#1}{\tilde}]{Y}_{\mkern -3mu (k)}}
\NewDocumentCommand{\Ytk}{s}{(\Tn[\IfBooleanT{#1}{\tilde}]{Y}_t)_{\mkern -3mu (k)}}
\NewDocumentCommand{\Yhk}{s}{(\Yh)_{(k)}}
\NewDocumentCommand{\Yhtk}{s}{(\Yht)_{(k)}}
\NewDocumentCommand{\Yv}{s}{\Vc[\IfBooleanT{#1}{\tilde}]{Y}}

\NewDocumentCommand{\KT} { s } {
  \llbracket 
  \IfBooleanTF{#1}{\lvec;}{}
  \Ak{1}, \dots,  \Ak{d} \rrbracket
}

\NewDocumentCommand{\KTt} { G{\st} } {
  \llbracket 
  #1;
  \Ak{1}, \dots,  \Ak{d} \rrbracket
}
\NewDocumentCommand{\KTtt} { G{\st}G{t} } {
  \llbracket 
  #1;
  \Akt{1}{#2}, \dots,  \Akt{d}{#2} \rrbracket
}

\NewDocumentCommand{\KThold} { } {
  \llbracket 
  \sh;
  \Akold{1}, \dots,  \Akold{d} \rrbracket
}

\NewDocumentCommand{\KTdef}{O{} O{} m !g}{%
  #1\llbracket \IfValueTF{#4}{#3;#4}{#3} #2\rrbracket%
}
\NewDocumentCommand{\KTauto}{m !g}{\KT[\left][\right]{#1}{#2}}
\NewDocumentCommand{\KTbig}{m !g}{\KT[\bigl][\bigr]{#1}{#2}}
\NewDocumentCommand{\KTbigg}{m !g}{\KT[\biggl][\biggr]{#1}{#2}}
\NewDocumentCommand{\KTBig}{m !g}{\KT[\Bigl][\Bigr]{#1}{#2}}
\NewDocumentCommand{\KTBigg}{m !g}{\KT[\Biggl][\Biggr]{#1}{#2}}

\NewDocumentCommand{\Gk} { G{k} } {%
  \Mx{G}^{(#1)}
}

\NewDocumentCommand{\Zk} { G{k} t'} {\MxRoot{}{Z}{#2}_{#1}}

\NewDocumentCommand{\Zkold} { G{k} t'} {\MxRoot{\lastiter}{Z}{#2}_{#1}}

\NewDocumentCommand{\Z} { t' } {\MxRoot{}{Z}{#1}}

\NewDocumentCommand{\I} {} {\mathcal{I}}

\NewDocumentCommand{\FPD} { s m m } {
  \IfBooleanTF{#1}
  {\partial #2 / \partial #3}
  {\frac{\partial #2}{\partial #3}}
}

\DeclareMathOperator{\diag}{diag}
\NewDocumentCommand \tol {} {\text{tol}}
\NewDocumentCommand \fest {} {\hat F}
\NewDocumentCommand \festold {} {\hat F_{\text{old}}}

\graphicspath{{./data/}}

\numberwithin{theorem}{section}

\newcommand{\TheTitle}{Streaming Generalized Canonical Polyadic Tensor Decompositions} 

\newcommand{\TheAuthors}{E. Phipps, N. Johnson and T. Kolda}

\headers{Streaming GCP Tensor Decompositions}{\TheAuthors}

\title{{\TheTitle}\thanks{%Submitted to the editors DATE.
\funding{Sandia National Laboratories is a multimission laboratory managed and operated by National Technology and Engineering Solutions of Sandia LLC, a wholly owned subsidiary of Honeywell International Inc., for the U.S. Department of Energy's National Nuclear Security Administration under contract DE-NA0003525.  This paper describes objective technical results and analysis. Any subjective views or opinions that might be expressed in the paper do not necessarily represent the views of the U.S. Department of Energy or the United States Government.  SAND2021-13340 R}}}

\author{
  Eric Phipps\thanks{Sandia National Laboratories, Albuquerque, NM
    (\email{etphipp@sandia.gov, nicjohn@sandia.gov}).}
  \and
  Nick Johnson\footnotemark[2]
  \and
  Tamara G. Kolda\thanks{Sandia National Laboratories, Livermore, CA
    (\email{tgkolda@sandia.gov}).}
}

\ifpdf
\hypersetup{
  pdftitle={\TheTitle},
  pdfauthor={\TheAuthors}
}
\fi

\begin{document}

\maketitle

% REQUIRED
\begin{abstract}
% !TEX root = online_gcp_paper.tex
In this paper, we develop a method which we call OnlineGCP for computing the Generalized Canonical Polyadic (GCP) tensor decomposition of streaming data.  GCP differs from traditional canonical polyadic (CP) tensor decompositions as it allows for arbitrary objective functions which the CP model attempts to minimize.  This approach can provide better fits and more interpretable models when the observed tensor data is strongly non-Gaussian.  In the streaming case, tensor data is gradually observed over time and the algorithm must incrementally update a GCP factorization with limited access to prior data.  In this work, we extend the GCP formalism to the streaming context by deriving a GCP optimization problem to be solved as new tensor data is observed, formulate a tunable history term to balance reconstruction of recently observed data with data observed in the past, develop a scalable solution strategy based on segregated solves using stochastic gradient descent methods, describe a software implementation that provides performance and portability to contemporary CPU and GPU architectures and integrates with Matlab for enhanced useability, and demonstrate the utility and performance of the approach and software on several synthetic and real tensor data sets.
\end{abstract}

% REQUIRED

% REQUIRED
%\begin{AMS}
%  %68Q25, 68R10, 68U05
%\end{AMS}

%\section*{Outline}
%\input{outline}

% The old writeup used as a starting point
%THIS IS THE OLD WRITEUP STARTING HERE:
%\input{old_writeup/writeup.tex}

\section{Introduction}
\label{sec:intro}
% !TEX root = online_gcp_paper.tex

We consider the problem of computing a generalized canonical polyadic (GCP) tensor decomposition \cite{HoKoDu20,KoHo20} in the situation where data is \emph{streaming}.
%The idea of streaming analysis is that you cannot see the entire dataset at once, only a portion.
% This may be because new data arrives incrementally.
% In other scenarios, you can take multiple passes over the data but it is too large to all fit in memory.
Generally speaking, the streaming paradigm assumes algorithms must update with a limited amount of data and limited passes on that data.
%Streaming tends to be an overloaded term that encompasses different sets of assumptions beyond data access constraints.
%
The data may be streaming because the volume of data is too large to fit in memory all at once;
however, the more general case is that the data is streaming because it is temporal and so arrives incrementally.
For example, consider a tensor that captures crime statistics in the city of Chicago so that entry $(i,j,k)$ is
the number of crimes of type $i$, in neighborhood $j$, at hour $k$.
In the streaming scenario, we receive a new 3-way tensor of crime statistics every day,
and we need to incorporate that information into the model.

The arrival of new data can be thought of in two different ways.
% tgk: Commenting out the reference to day mode since the tensor was just described without a day mode, so the reader will be confused.
% In one view, there is no ``day'' mode. 
We could view the Chicago Crime data as a 3-way tensor (type $\times$ neighborhood $\times$ hour)
with new observations each day that can be considered a \emph{statistical sample}.
Alternatively, we can have an explicit time mode.
For Chicago crime, the tensor is then a 4-way tensor
(type $\times$ neighborhood $\times$ hour $\times$ day) with a new hyperslice \emph{appended} daily.
In this latter case, the fourth mode corresponding to the day is growing and referred to as the \emph{temporal} mode.  The GCP tensor decomposition computes a \emph{factor matrix} for each mode,
so in the Chicago example we have a crime-type factor matrix for mode one, a neighborhood factor matrix for mode two, and an hour factor matrix for mode three, whether we think of it as a 3-way or 4-way tensor.
In the 4-way interpretation, we additionally have a day factor matrix for mode four,
and a new row is added to that factor matrix with each new day of data.
In the 3-way interpretation, we can think instead of adjusting the weights of the factors each day.
Ultimately, these two viewpoints are not very different since the new row in the temporal factor matrix in the 4-way interpretation is roughly equivalent to the updated weights in the 3-way interpretation.

A more interesting assumption is whether or not the underlying generative processes are changing with time.
There is always some balancing of new information and old.
If these processes are unchanging, then we may expect that our estimates of the non-temporal factor matrices will converge after a suitable number
of observations.
Such a formulation assumes data observed comes from some consistent, but unknown distribution.
This aligns better with incremental algorithms that progressively converge to a single value, and the ordering of observations aligns with sampling assumptions.
In this situation, it is often useful to give a heavy weight to older information, slowing the amount of change allowed as more observations accumulate. In fact, in these cases, the order that the information arrives is irrelevant.

In most cases of interest, however, the generative processes are changing as well, and we are interested in understanding these shifts.
This is sometimes referred to as \emph{concept drift}, where it is assumed the distributions of observed data are evolving in time \cite{gama2014survey}.
Algorithms must be designed to adapt with drifting data, and the ordering of observations is critically important to track the evolution of the data distribution.
In this case, we need to balance between adapting to changing generative processes without confusing them for statistical fluctuations.

Much research has been done in the case of the canonical polyadic (CP), also known as CANDECOMP/PARAFAC, decomposition.
In this work, we extend existing methods to GCP, which differs from CP in that GCP allows for
arbitrary objective functions.  In particular, our contributions are as follows:
\begin{itemize}
	\item As there is no single streaming problem formulation in the literature, we provide a concise overview of the streaming literation emphasizing the various assumptions made in different works (\cref{sec:related});
	\item We extend the GCP formalism to the streaming context by deriving a GCP optimization problem to be solved as each tensor slice is observed enabling CP factorization of tensors using arbitrary objective functions (\cref{sec:streaming_gcp_formulation});
	\item This formulation incorporates a tunable history term into the optimization problem to balance reconstruction of recently observed data with data observed in the past;
	\item We develop a solution strategy for the GCP streaming problem based on segregated solves of the temporal weights and factor matrices using stochastic gradient descent solution methods (\cref{sec:streaming_gcp_solution,sec:streaming_gcp_alg});
	\item We provide a highly performant software implementation of the algorithm in our GenTen software package through an extensible Matlab class hierarchy leveraging low-level math kernels implemented on top of the Kokkos framework providing scalable thread parallelism and portability to contemporary CPU and GPU architectures (\cref{sec:software});
	\item We demonstrate the utility and performance of the approach on a variety of synthetic and realistic data sets using several GCP objective functions (\cref{sec:numerical}).
\end{itemize}

%\tknote{Eric --- add bulleted list of contributions of this work here.}

%\subsection{Streaming}

%%% Local Variables:
%%% mode: latex
%%% TeX-master: "online_gcp_paper"
%%% End:

\section{Background and Notation}
\label{sec:background}
% !TEX root = online_gcp_paper.tex

In this work, we assume the reader to be generally familiar with tensors and tensor decomposition methods.
For a thorough overview, we refer the reader to Kolda and Bader~\cite{KoBa09}.
Following standard practice, we denote tensors by bold calligraphic letters (e.g., $\X$),
matrices by bold capital letters ($\Mx{A}$),
vectors by bold lowercase letters ($\Vc{a}$)
and scalars by lowercase letters ($a$).
We use multi-index notation to indicate tensor elements, i.e.,
$x_i \equiv \xi*$ denotes the entry $i=(i_1,\dots,i_d)\in \I \equiv \set{1, \dots, I_1} \otimes \dots \otimes \set{1, \dots, I_d}$ of the $d$-way tensor $\X\in\Real^{I_1 \times \cdots \times I_d}$.

\subsection{Canonical Polyadic (CP) Tensor Decompositions}
\label{sec:cp}

For a given $d$-way tensor $\X \in \Real^{I_1 \times \cdots \times I_d}$, the Canonical Polyadic (CP) decomposition,
also known as the CANDACOMP/PARAFAC decomposition,
attempts to find a good approximating low-rank model tensor $\M$ of the form 
\begin{equation}\label{eq:CP_model}
  \X \approx \M  = \sum_{j=1}^R \lj"  \akj[1] \circ \akj[2] \circ \dots \circ \akj[d]
\end{equation}
where $\lj$ is a scalar weight,
$\akj$ is a column vector of size $I_k$,
$\circ$ represents the tensor outer product, and $R$ is the approximate rank.
The column vectors for each mode $k$ are often collected into a matrix
$\Ak = [ \akj{1} \, \cdots  \; \akj{R} ]$
of size $I_k\times R$ called a factor matrix. 
Given a weight vector $\lvec = [ \lj{1} \, \cdots \, \lj{R} ]^T$ and factor matrices $\set{ \Ak{1}, \dots, \Ak{d} }$, we refer to the resulting low-rank model $\M$ as a Kruskal tensor (or K-tensor for short) and use the short-hand notation $\M=\KT*$~\cite{BaKo07}. 
For traditional CP decompositions, $\M$ is computed by solving a nonlinear least-squares problem
\begin{equation}\label{eq:CP_problem}
  \min_{\M} \;\; \| \X - \M \|_F^2 \quad \mbox{s.t.} \quad \M = \KT*
\end{equation}
where $\| \X - \M\|_F^2 = \sum_{i \in \I} (x_i - m_i)^2$, with $\I$ defined as above, denotes the tensor Frobenius (sum-of-squares) norm.
Note that in~\cref{eq:CP_problem}, the minimization is with respect to both the weights $\lvec$ and factor matrices $\set{ \Ak{1}, \dots, \Ak{d} }$.
Many approaches have been developed for efficiently solving \cref{eq:CP_problem} that are scalable to large, sparse tensors. 
However, a very common, successful approach that is also relevant to the streaming problem is alternating least-squares (ALS) which is an iterative method, that for each iteration, cycles over modes $k=1,\dots,d$, holds all of the modes other than mode $k$ fixed, and solves the resulting linear least squares problem for $\Ak$.

\subsection{Generalized CP Decompositions}
\label{sec:gcp}

As described in~\cite{HongKoldaDuersch2020}, the CP problem~\eqref{eq:CP_problem} is equivalent to a maximum likelihood estimation procedure where the entries $x_i$ of the tensor of $\X$ are i.i.d. 
Gaussian with with mean $m_i$ and some variance $\sigma^2$ which is constant across the tensor, i.e., $x_i = \mathcal{N}(m_i,\sigma)$. 
Such a statistical assumption may not be appropriate for many types of data (e.g., count or binary),  motivating the development of the Generalized Canonical Polyadic (GCP) method~\cite{HongKoldaDuersch2020}. 
 In this method, it is assumed the tensor entries follow some known, parameterized probability distribution $x_i \sim p(x_i|\eta_i)$ determining the likelihood of each entry $x_i$, where $\eta_i$ is the (unknown) parameter of the distribution.  In this case, the CP model is computed to maximize the likelihood $p(x_i|\eta_i)$ of the tensor entry observation $x_i$ through an invertible link function $\ell(\eta_i) = m_i$ connecting the CP model parameter $m_i$ to the distributional parameter $\eta_i$.  This results in the more general optimization problem
\begin{equation}\label{eq:GCP_problem}
  \min_{\M} \;\; F(\X,\M) = \sum_{i \in \I} f(x_i,m_i) \quad \mbox{s.t.} \quad \M = \KT*,
\end{equation}
where as before $\I = \set{1, \dots, I_1} \otimes \dots \otimes \set{1, \dots, I_d}$ is the set of all tensor multi-indices (including both zeros and nonzeros).  Here $f(x,m) = -\log p(x|\ell^{-1}(m))$ is the negative log-likelihood and is called the loss function.  For example, one may have $f(x,m) =  m - x\log m$ with $\ell(\eta) = \eta$ for a tensor containing count data, assuming a Poisson distribution, or $f(x,m) = \log(m+1) - x\log m$ with $\ell(\eta) = \eta/(1-\eta)$ for a binary tensor assuming a Bernoulli distribution.\footnote{In practice, $\log m$ is replaced by $\log(m+\epsilon)$ where $\epsilon$ is a small constant to allow $m=0$.  Also, depending on the choice of loss function, the minimization problem~\cref{eq:GCP_problem} may include additional constraints such as $m_i\geq0$.}  See~\cite{HongKoldaDuersch2020} for a detailed derivation of these loss functions for different statistical distributions. 
In the Gaussian case, $f(x,m) = (x-m)^2$, so~\cref{eq:CP_problem} becomes a special case.  It is important to note that in the general case however, the CP model no longer approximates the tensor itself, but rather the natural parameter of the distribution underlying the assumed statistical model of the tensor data, which will be crucial for the streaming method described later.

The challenge in the GCP method is solving~\cref{eq:GCP_problem} for general loss functions $f$ which loses the least-squares structure, making ALS-type approaches impossible. 
In~\cite{HongKoldaDuersch2020}, the authors instead pursue gradient-based optimization approaches and derive the corresponding gradient formulas:
\begin{align}
  \FPD{F}{\Ak} &= \Yk \Zk \diag(\lvec), \quad k=1,\dots,d, \label{eq:gcp_grad_A}\\
  \FPD{F}{\lvec} &= \Z' \Yv \label{eq:gcp_grad_lambda}
\end{align}
where $\Y\in\Real^{I_1\times\dots\times I_d}$ is a gradient tensor defined by
\begin{equation}
  y_i = \frac{\partial f}{\partial m}(x_i,m_i), \quad i\in\I.
\end{equation}
Here $\Yk$ denotes the mode-$k$ matricization/unfolding of $\Y$, $\Yv=\mbox{vec}(\Y)$ is the vectorization of $\Y$,
\begin{align}
  \Zk &= \Ak{d}\odot\dots\odot\Ak{k+1}\odot\Ak{k-1}\odot\dots\odot\Ak{1}, \quad k=1,\dots,d, \\
  \Z &= \Ak{d}\odot\Ak{d-1}\odot\dots\odot\Ak{1},
\end{align}
and $\odot$ denotes the Khatri-Rao product. 
Thus the factor matrix gradients~\eqref{eq:gcp_grad_A} are given by the Matricized Tensor Times Khatri-Rao Product (MTTKRP) involving the gradient tensor $\Y$. 
Note that $\Y$ is in general dense, even if $\X$ is sparse, making traditional gradient-based methods impractical for large, sparse $\X$. 
Instead, the authors in~\cite{KoHo20} leverage stochastic gradient descent (SGD) in this case, using randomly sampled gradients of the form 
\begin{align}
  \FPD{F}{\Ak} &\approx \Yk* \Zk \diag(\lvec), \quad k=1,\dots,d, \\
  \FPD{F}{\lvec} &\approx \Z' \Yv*
\end{align}
where $\Y*$ is a sparse, randomly sampled approximation of $\Y$. 
In the sequel, we will leverage these formulas for developing the streaming GCP algorithm for sparse tensors, employing the sparse, stratified sampling methodology of~\cite{KoHo20}.

%%% Local Variables:
%%% mode: latex
%%% TeX-master: "online_gcp_paper"
%%% End:

\section{Related Work}
\label{sec:related}
% !TEX root = online_gcp_paper.tex

We review the work in the domain of streaming or online CP tensor decomposition.
There is no single well-defined problem in this context, so we try to explain
the different formulations and assumptions.
We make a few assumptions throughout.
\begin{itemize}
\item  Updates are processed in \emph{discrete} batches indexed by time $t=1,2,\dots$.
\item At each time $t$, a $d$-way tensor is observed, possibly incomplete.
\item Dimensions are fixed throughout all time, unless otherwise stated.
\item The CP rank is known and fixed, unless otherwise stated.
\end{itemize}

\subsection{Problem setup for two-way temporal slices}
\label{sec:problem-setup-d=2}

In the case that $d=2$, we receive a matrix $\Xt* \in \Real^{I \times J}$
for each time $t=1,2,\dots$.

If the temporal mode is finite so that $t=1,\dots, T$, we can consider
that $\X$ is the $I \times J \times T$ tensor formed by stacking all time slices.
For a given rank $R$, the standard goal is to find factor matrices
$\A \in \Real^{I \times R}$, $\B \in \Real^{J \times R}$, and $\Mx{S} \in \Real^{T \times R}$
that minimize
\begin{displaymath}
  \min_{\A,\B,\Mx{S}} \sum_{i=1}^I \sum_{j=1}^J \sum_{t=1}^T (x_{ijt} - m_{ijt})^2
  \qtext{s.t.}
  m_{ijt} = \sum_{\ell=1}^R a_{i\ell} b_{j\ell} s_{t\ell}.
\end{displaymath}
There are a few different streaming and online formulations of this problem.

In one formulation, the challenge is that we can only see one (or a few) temporal slices $\Xt*$ at any given time.
This may be due to memory constraints.
However, we assume that the factor matrices $\A$ and $\B$ are fixed for all time.

In other versions, the factor matrices $\A$ and $\B$ can change over time.
Then the problem becomes more interesting. 
Let $\st \in \Real^R$ denote row $t$ (transposed) of $\Mx{S}$, i.e., 
\begin{equation}\label{eq:S-matrix}
  \Mx{S} =
  \begin{bmatrix}
    \st{1} & \cdots & \st{T}
  \end{bmatrix}^{\intercal}.
\end{equation}
Then, at time step $t$, the goal is to find $\st \in \Real^R, \A \in \Real^{I \times R}, \B \in \Real^{J \times R}$ that minimize
\begin{displaymath}
 \min_{\A,\B,\st} \|\Xt* - \Mt*\|^2 \equiv \sum_{i=1}^I \sum_{j=1}^J (x_{ijt} - m_{ijt})^2
  \qtext{s.t.}
  \Mt* = \A \diag(\st)\B'
\end{displaymath}
If we only fit $\Xt*$, however, the problem is not well defined, i.e.,
it does not produce essentially unique minimizers $\A$ and $\B$.
Instead, at time $t$, it is common to include some \emph{historical} information in the objective function,
the exact details of which depend on the formulation.

It can be argued that a truly streaming problem is not finite, so $t=1,2,\dots$.
In that case, we cannot save all the historical information.
Additionally, such problems are generally more interesting if the factor matrices change slowly in time;
otherwise, we can assume that the factor matrices would be learned within finite time and the
only thing changing at each time step are the weights $\st$.

\subsection{Problem setup for higher-order temporal slices}
\label{sec:problem-setup-d2}

If $d > 2$, the updates are tensors. For each time $t=1,2,\dots$,
we receive a tensor $\Xt \in \Real^{I_1 \times I_2 \times \cdots \times I_d}$.
At time $t$, 
the goal is to find factor matrices $\Ak \in \Real^{I_k \times R}$ for $k=1,\dots,d$ and
weights $\st \in \Real^R$ that minimize
\begin{equation}\label{eq:single-time-step}
 \min_{\st,\Ak{1},\dots,\Ak{d}} \| \Xt-\Mt\|^2 \equiv \sum_{i\in\I} (\xit - \mit)^2
  \qtext{s.t.}
  \Mt = \KTdef{\st}{\Ak{1},\dots,\Ak{d}},
\end{equation}
generally with some methodology for incorporating historical information.
Here we use the shorthand
%\begin{equation}
%  \label{eq:I-def}
%  \I \equiv \set{1,\dots,I_1} \otimes \set{1,\dots,I_2} \otimes  \cdots \otimes \set{1,\dots,I_d},
%\end{equation}
%and 
$\xit \equiv \Xt(i_1,\dots,i_d)$ and $\mit \equiv \Mt(i_1,\dots,i_d)$.

\subsection{Earliest work}
\label{sec:earliest-work}

To the best of our knowledge, the earliest work in this area is Nion and Sidiropoulos \cite{NiSi09}.
They consider the case where each observation is a two-way matrix, as in \cref{sec:problem-setup-d=2}.
The rank and sizes are fixed across time, and the factors are assumed to be \emph{slowly} varying, though
there are no experiments with factors that varied in time in that work.
Their primary focus was on demonstrating %that this was 
an alternative method for fitting CP decompositions
that traded a small amount of accuracy for increased speed.
Their general formulation of the problem is as follows.
At time $t$, solve
\begin{equation}\label{eq:NiSi09}
  \min_{\A,\B,\st} \sum_{h=1}^{t} \theta^{t-h} \| \Xh* - \Mh* \|^2
  \qtext{s.t.} \Mh* = \A'\!\diag(\sh) \B \text{ for all } h \in \set{1,\dots,t}.
\end{equation}
The parameter $\theta \in (0,1)$ downweights older hyperslices.
The weights $\sh$ are fixed for all $h<t$; however, the $\A$ and $\B$ matrices
are updated at each time step.
This means that the $\A$ and $\B$ should still somewhat fit the older data.
The summation over time incorporates \emph{historical} information.
We omit the details of the method since it is relatively complicated and
has been subsequently bested by other methods.
Technically, this method requires all historical information.
However, because $\theta^{t-h}$ is exponentially decreasing, older information
can effectively be discarded after a small number of time steps.

\subsection{Online SGD}
\label{sec:online-sgd}

Mardani, Mateos, and Giannakis \cite{MaMaGi15} consider both matrix and 3-way tensor streaming;
we discuss only the tensor streaming part of their work.
In contrast to \cref{eq:NiSi09}, Mardani et al.~\cite{MaMaGi15} account for missing data, add regularization, and
have an entirely different computational approach.

To account for missing data, define the matrix
\begin{displaymath}
  \Wt*(i,j) =
  \begin{cases}
    1 & \text{if entry $(i,j)$ is known at time $t$},\\
    0 & \text{otherwise}.
  \end{cases}
\end{displaymath}
Additionally, Mardani et al.\@ add regularization with parameter $\lambda$.
At time $t$, the formulation is
\begin{multline}\label{eq:MaMaGi15}
  \min_{\A,\B,\st} \sum_{h=1}^{t } \theta^{t-h} \| \Wh*\ast( \Xh* - \Mh* ) \|
  + \bar \lambda_t (\|\A\|_F^2 + \|\B\|_F^2) + \lambda \|\st\|_2^2\\
  \qtext{s.t.} \Mh* = \A'\!\diag(\sh) \B \text{ for all } h \in \set{1,\dots,t}.
\end{multline}
As with \cref{eq:NiSi09}, the parameter $\theta \in (0,1)$ downweights older hyperslices.
For writing efficiency, we pulled the term for time $t$ into the summation, but
the weights $\sh$ are fixed for all $h<t$.
The asterisk ($\ast$) denotes \emph{elementwise} multiplication,
and the effect is that only observed entries are included in the summation.
The definition of $\bar \lambda_t$ is somewhat unclear in the paper. %due to many typos.
At one point,
it seems to propose that $\bar\lambda_t \equiv \lambda/\sum_{h=1}^t \theta^{t-h}$ to
ensure a degree of consistent weighting as compared to the regular 3-way problem in the finite case,
but the pseudo-code seems to indicate in the code that either
$\bar\lambda_t \equiv \lambda/t$ or $\bar\lambda_t \equiv \lambda/(t\sum_{h=1}^t \theta^{t-h})$.
In the experiments, the regularization parameter is set according
to a standard in matrix completion: $\lambda = \sqrt{2IJ\pi}\sigma$ where $\pi$ is the proportion of sampled data at
each time step and $\sigma$ is the noise level.
%\epnote{If you look at equation P5 of the OnlineSGD paper, it looks like the sum over time slices in \cref{eq:MaMaGi15} includes the two regularization terms, which are also multiplied by the downweighting factor, and the last term should be $\lambda\|\sh\|_2^2$, i.e., it should read $\min_{\A,\B,\st} \sum_{h=1}^{t } \theta^{t-h} \left[\| \Wh*\ast( \Xh* - \Mh* ) \| + \bar \lambda_t (\|\A\|_F^2 + \|\B\|_F^2) + \lambda \|\sh\|_2^2\right]$.  In this case, the definition of $\bar\lambda_t = \lambda/\sum_{h=1}^t \theta^{t-h}$ makes sense, because it is then equivalent to moving that term out of the sum.  Also, as a minor point, they also use $\lambda_t$ instead of $\lambda$, to allow for time-varying regularization, but I don't think they actually do that anywhere.  But then if you look at Eq. 18, they appear to have moved the factor matrix regularization terms out of the sum (without changing the definition of $\bar\lambda_t$) and dropped the spatial mode factorization term.  Then finally in Algorithm 3, they drop the sum altogether and just consider the time $t$ optimization.}\tkwarning{Agreed, but later they put a $t$ in the denominator, starting at the last line on page 2670 and also in Alg.~3.}

Problem \eqref{eq:MaMaGi15} is solved iteratively with two basic steps at each iteration.
First, $\st$ is solved for via a closed form expression holding $\A$ and $\B$ fixed.
Second, the method takes \emph{one step} of gradient descent (GD) for updating $\A$ and $\B$.
They call this \emph{stochastic} gradient descent (SGD) because there only partial data at each step.

There are some implicit assumptions that are not clearly stated in the paper.
There is no proof that the stochastic gradient is correct in expectation.
To do so requires some assumptions about how the data is sampled
and also appropriate weighting.
The experiments (on cardiac dynamic MRI and Internet traffic) are constructed so that it all works correctly enough --- the data is sampled
uniformly and the same number of samples are taken at each time step.

\subsection{CP-Stream}
\label{sec:cp-stream}

CP-Stream \cite{SmHuSiKa18} is similar to OnlineSGD \cite{MaMaGi15}.
The primary difference is that it avoids saving the older data and
instead uses an approximation.
It also works for $d>2$.
At time step $t$, CP-Stream executes two phases.
Let $\Akold$ denote the \emph{old} factor matrices, i.e., from time $t-1$.
The first phase computes $\st$ with all the old factor
matrices by solving
\begin{equation}
  \label{eq:cp-stream-step1}
  \min_{\st} \| \Xt - \Mtold \|^2 + \lambda \|\st\|^2 \qtext{s.t.} \Mtold \equiv  \KTdef{\st}{\Akold{1},\dots,\Akold{d}}.
\end{equation}
The second phases computes $\set{\Ak{1},\dots,\Ak{d}}$,  estimating the observations from prior time steps via $\Xh \approx \Mhold \equiv \KTdef{\sh}{\Akold{1},\dots,\Akold{d}}$:
\begin{multline}
  \label{eq:cp-stream-step2}
  \min_{\Ak{1},\dots,\Ak{d}} \| \Xt - \Mt \|^2 + \sum_{h=1}^{t-1} \theta^{t-h} \| \Mhold - \Mh \|^2 \\
  \qtext{s.t.} \Mh =  \KTdef{\sh}{\Ak{1},\dots,\Ak{d}} \text{ for all } h=1,\dots,t.
\end{multline}
%\tknote{Need to double-check the regularization, handling missing data, etc.}

\subsection{Other works}
\label{sec:other-works}

%\tknote{I haven't mentioned that GridTF \cite{PhCi11a} approach, but feel free to add it in.}

\paragraph{OnlineCP}
%\label{sec:onlinecp}

The OnlineCP method \cite{ZhViBaJi16} works for arbitrary $d$-way tensors
and solves exactly the standard least squares subproblems,
regularizing by taking only a single step of ALS.
%TGK: I changed it! Thanks for catching the error.
%\epnote{I don't think this is exactly correct.  OnlineCP does include prior slices using the clever approach mentioned in the next sentence.}
The innovation is the clever reuse of expensive calculations when folding
in each time slice. It effectively assumes that the factor matrices are fixed.

\paragraph{OLSTEC}
%\label{sec:olstec}

The Online Low-Rank Subspace Tracking by Tensor CP Decomposition
(OLSTEC) method \cite{Ka16} is similar to one of the methods proposed
in \cite{NiSi09}, but it can handle missing data and is the first paper to consider
changes in the factor matrices in its experimental results.
The experiment results show that they
do better in this regime than OnlineSGD \cite{MaMaGi15}.

\paragraph{MAST}

Multi-aspect Streaming Tensor (MAST) \cite{SoHuGeCa17} is
notable because it allows for the non-temporal modes to grow in
time. (It considers both CP and Tucker.)

\paragraph{SamBaTen}

Sampling-Based Incremental Tensor Decomposition (SamBaTen) \cite{GuPaPa18} 
samples multiple subtensors, factors those independently, and then merges the results.
It depends heavily on the results being essentially unique and consistent across the subtensors,
which necessarily assumes that the factors are not changing in time.
The temporal aspect is not clear since the entire tensor (across all time) seems to be saved.

\paragraph{SeekAndDestroy}

SeekAndDestory~\cite{PaGuPa19} handles concept drift by allowing the addition of new factors
as time progresses. It is not specifically a streaming algorithm because it is not updating
the factorization so much as augmenting it.
It receives data in batches, computes the CP decomposition from scratch,
and then it merges this with the information from prior batches.
We do not consider this to be a streaming method because the existing decomposition is not
updated directly.
Instead, SeekAndDestroy finds those factors that
are overlapping and then identifies and older factors that do not appear in the new batch as well as any
new factors in the new batch. For rank determination on each batch, it uses a heuristic called AutoTen.
This method depends heavily on each new batch having sufficient information to compute a full and
essentially unique decomposition as this is the only way to ensure that overlap with past factors
can be identified.

The ENSIGN software \cite{LeBaHeEz18} implements a method similar to SeekAndDestroy.
In addition to CP-ALS, they include CP-APR and CP-ALS-NN.

\paragraph{Bayesian Methods}

Probabilistic Stream Tensors (POST) \cite{DuZhLeZh18} and
Variational Bayesian Inference (VBI) \cite{ZhHa18} are two papers
that propose priors for the tensor model. POST considers models
for both continuous and binary data.
VBI models each time step as a CP model plus sparse noise ($\Tn{S}_t$)
and Gaussian noise ($\Tn{E}_t$):
\begin{gather*}
  \label{eq:vbi}
  \Xt = \KTdef{\st; \Ak{1},\dots,\Ak{d}} + \Tn{S}_t + \Tn{E}_t, \\
  \begin{aligned}
  \Ak(i_k,j) &\sim \mathcal{N}(0,\lambda_j) \\
  \st(j) &\sim \mathcal{N}(0,\lambda_j) &
  \lambda_j &\sim \text{InvGamma}(\alpha_{\lambda},\beta_{\lambda})\\
  \Tn{S}_t(i_1,i_2,\dots,i_d) &\sim \mathcal{N}(0,\gamma_{i_1i_2\cdots i_d}) &
  \gamma_{i_1i_2\cdots i_d} & \sim \text{InvGamma}(\alpha_{\gamma},\beta_{\gamma})\\
  \Tn{E}_t(i_1,i_2,\dots,i_d) &\sim \mathcal{N}(0,\tau) &
  \tau & \sim \text{InvGamma}(\alpha_{\tau},\beta_{\tau})    
  \end{aligned}
\end{gather*}
The prior on the $\lambda$'s encourages low-rank, and prior on the $\gamma$'s allows for
sparse outliers. The method is amenable to missing data as well.
Both methods are very slow to compute and are not scalable.

%%% Local Variables:
%%% mode: latex
%%% TeX-master: "online_gcp_paper"
%%% End:

\section{Streaming GCP}
\label{sec:streaming}
% !TEX root = online_gcp_paper.tex

We now consider GCP factorization in the streaming context.
We first motivate and describe the minimization problem for the streaming GCP problem, then describe the solution strategy, and then conclude with a summary of the solution algorithm which we call OnlineGCP.

\subsection{Streaming GCP problem formulation}
\label{sec:streaming_gcp_formulation}

We are primarily interested in the infinite streaming problem where at each time step $t$, a new $d$-dimensional tensor $\Xt\in\Real^{I_1\times\dots\times I_d}$ is observed.
We assume the dimensions $I_1,\dots,I_d$ of each tensor do not change.
Since both the non-streaming GCP and OnlineSGD solution algorithms rely on gradient descent, our approach for streaming GCP is inspired by OnlineSGD.

For simplicity, we begin with the problem for a single temporal slice, $\Xt$,
and add on history in the discussion that follows.
For just time slice $t$, we pose the optimization problem for given rank $R$ as
\begin{equation}\label{eq:gcp_stream_t}
\begin{gathered}
  \min_{\Mt} \;\;
  \sum_{i \in \I} f(\xit,\mit) +
  \frac{\lambda}{2} \sum_{k=1}^d \| \Ak \|_F^2 + \frac{\mu}{2} \| \st \|_2^2 \\
  \mbox{s.t.} \quad \Mt = \KTt, \;\; \Ak{1},\dots,\Ak{d},\st \ge l.
\end{gathered}
\end{equation}
As above, we use the shorthand
$\xit \equiv \Xt(i_1,\dots,i_d)$ and $\mit \equiv \Mt(i_1,\dots,i_d)$,
$\Ak \in \Real^{I_k \times R}$ for $k = 1,\dots,d$ are the factor matrices, and $\st \in \Real^{R}$ is the weight for time $t$.
We incorporate the option for a lower bound $l$ on the factor matrix/weight entries
(with the understanding that $l=0$ for nonnegativity constraints and $l=-\infty$ for problems where there is no lower bound).
As in OnlineSGD, we include regularization terms for the factor matrices $\Ak$ and weights $\st$
with multipliers $\lambda$ and $\mu$, respectively, to encourage low-rank~\cite{6579771}.
We do not explicitly include the dependency of the factor matrices on $t$ since they are, ideally, less
sensitive to time. 

In general, we want to incorporate historical information to keep the factor matrices from changing too much
at each time step.
There are many ways such historical information could be included, and several of these have
been used in the previous work discussed in \cref{sec:other-works}.
One possible approach is to add \emph{historical regularization} to \cref{eq:gcp_stream_t} via the second term in the following
where we define the historical model to be the old weight with the current factor matrices, i.e.,
$\Mh \equiv \KTt{\sh}$:
\begin{equation}\label{eq:gcp_stream_problem_1}
\begin{gathered}
  \min_{\Mt} \;\; \sum_{i \in \I} f(\xit,\mit)
  + \sum_{h\in\HWt} w_h \sum_{i\in\I} f(\xih,\miht)
  + \frac{\lambda}{2} \sum_{k=1}^d \| \Ak \|_F^2 + \frac{\mu}{2} \| \st \|_2^2 \\
  \mbox{s.t.} \quad \Mt = \KTt, \;\;
  \Ak{1},\dots,\Ak{d},\st \ge l.
\end{gathered}
\end{equation}
Here, $\HWt \subseteq \{1,\dots,t-1\}$.
The terms $\xih=\Xh(i_1,\dots,i_d)$ and $\miht = \Mh(i_1,\dots,i_d)$
index the ``historical'' tensors.
The weights $w_h$  control the importance of historical terms.

The time index $t$ could be infinite, so we limit ourselves to a history window
of fixed size such that $|\HWt| = \min\set{t-1,H}$ for some fixed constant $H$.
Moreover, we impose $\HWt \subseteq \HWt{t-1} \cup \set{t-1}$
so that no older information is ever added to the history.
\emph{This means that we can discard all older information except what's
in the history.}
There are many ways in which the history window and weighting could be chosen.
In our work, we use 
reservoir sampling~\cite{VitterTOMS1985} which works as follows.
For $t\leq H+1$ we set $\HWt=\{1,\dots,t-1\}$.  Then for $t > H+1$, We set $\HWt=\HWt{t-1}$ with probability $1-H/(t-1)$;
otherwise, we set $\HWt$ to be $\HWt{t-1}$ where we have ejected one existing element and replaced it with $t-1$.
%\tknote{What precisely happened in the first $H$ time steps? Is $\HWt{H} = \set{1,\dots,t}$?}
This ensures that $\HWt$ is a uniform random sample of $\set{1,\dots,t-1}$.
We use exponential weighting of the form $w_h = w \theta^{t-h}$ where $0<\theta\leq 1$
and $w$ is a multiplier allowing the entire history term to be scaled by a constant.

The approach outlined so far stores at most $H$ temporal slices, which could require significant memory storage if a large window size $H$ is desired.
Following the approach of CP-Stream (see \cref{sec:cp-stream}),
we would like to reduce storage costs further in the case where the factor matrices are assumed to change slowly in time by approximating these slices using the factor matrices $\Akold{1},\dots,\Akold{d}$ from the previous time step.  It is tempting to replace $\xih$ in \cref{eq:gcp_stream_problem_1} with $\mihold$ where $\Mhold \equiv\KThold$ is the CP model derived from the prior factor matrices and the historical weights, and $\mihold = \Mhold(i_1,\dots,i_d)$.  However, this is not a valid approximation in general because, as described in \cref{sec:gcp}, the CP model generated in the GCP method does not directly approximate the data tensor, but instead the parameter of the assumed probability distribution (in fact, the support of $\xih$ may not even coincide with $\mihold$, e.g., for a binary tensor, $\xih\in\{0,1\}$ whereas $\mihold\in(0,\infty)$).  Instead, we propose adding a historical regularization term that penalizes changes in the CP model using the Frobenius norm:
\begin{equation}\label{eq:gcp_stream_problem_2}
\begin{gathered}
  \min_{\Mt} \;\; \sum_{i \in \I} f(\xit,\mit)
  + \frac{1}{2} \sum_{h\in\HWt} w_h \|\Mhold-\Mh\|_F^2
  + \frac{\lambda}{2} \sum_{k=1}^d \| \Ak \|_F^2 + \frac{\mu}{2} \| \st \|_2^2 \\
  \mbox{s.t.} \quad \Mt = \KTt, \;\; \Ak{1},\dots,\Ak{d},\st \ge l
\end{gathered}
\end{equation}
We reiterate that in \cref{eq:gcp_stream_problem_2}, the optimization is over the factor matrices $\Ak{1},\dots,\Ak{d}$ and temporal weights $\st$ for each time step $t$, with $\Akold{1},\dots,\Akold{d}$ and $\sh$ for $h\in\HWt$ held fixed.\footnote{We note that a third possible approach is to directly penalize changes in the factor matrices by replacing the history regularization term in \cref{eq:gcp_stream_problem_2} with a term of the form $\sum_{k=1}^d \tilde{w}_k \|\Akold-\Ak\|_F^2$.  However such an approach requires careful tuning of the regularization parameters $\tilde{w}_k$ since it does not incorporate the historical weights $\sh$.}

\subsection{Streaming GCP solution strategy}
\label{sec:streaming_gcp_solution}

We now describe the proposed solution strategy for~\cref{eq:gcp_stream_problem_2}, which we call OnlineGCP.
Assume~\cref{eq:gcp_stream_problem_2} has been solved for $h=1,2,\dots,t-1$ resulting in the current approximations $\Akold{1},\dots,\Akold{d}$, $\Vc{s}_{t-1}$, with $\sh$ for $h\in\HWt{t-1}$ already known from prior iterations.  First choose the new history window $\HWt\subseteq\HWt{t-1}\cup\set{t}$.  Given the new tensor slice $\Xt$, define
\begin{multline}\label{eq:gcp_stream_objective}
  F(\Xt,\Ak{1},\dots,\Ak{d},\st) = \\
  \sum_{i \in \I} f(\xit,\mit)
  + \frac{1}{2}\sum_{h\in\HWt} w_h \|\Mhold-\Mh\|_F^2
  + \frac{\lambda}{2} \sum_{k=1}^d \| \Ak \|_F^2 + \frac{\mu}{2} \| \st \|_2^2
\end{multline}
to be the streaming GCP objective function.
We then solve \cref{eq:gcp_stream_problem_2} using a two-step minimization procedure inspired by OnlineSGD.  In particular, we first solve \cref{eq:gcp_stream_problem_2} for $\st$ with $\Ak{k}=\Akold{k}$ held fixed, namely
\begin{equation}\label{eq:gcp_stream_temporal_problem}
\begin{gathered}
  \min_{\st} F(\Xt,\Mt) = 
  \sum_{i \in \I} f(\xit,\mit) + 
  %\frac{\lambda}{2} \sum_{k=1}^d \| \Akold{k} \|_F^2 + 
  \frac{\mu}{2} \| \st \|_2^2 \quad
  \mbox{s.t.}\quad \st \ge l.
\end{gathered}
\end{equation}
The history term and factor matrix regularization terms are dropped because they have no dependence on $\st$.
To solve \cref{eq:gcp_stream_temporal_problem}, we use the SGD solver described in \cref{sec:gcp} for the static GCP problem,
modified to only solve for $\st$ with the factor matrices fixed.
Leveraging~\cref{eq:gcp_grad_lambda}, the gradient for this subproblem is
\begin{equation}
  \FPD{F}{\st} = \Z' \Yv_t + \lambda \st,
\end{equation}
where as before $\Yv_t = \mbox{vec}(\Yt)$ and $\Yt \in \Real^{I_1 \times \dots \times I_d}$ is the gradient tensor for slice $\Xt$ defined by $y_{it} = \FPD{f}{m}(\xit,\mit)$.
As in the static case, this tensor is sampled each SGD iteration resulting in each stochastic gradient.

Once $\st$ is computed, the factor matrices $\Ak$ are computed %from $\Akold{k}$
by applying a fixed number of SGD iterations to \cref{eq:gcp_stream_objective}, holding $\st$ fixed.
Using~\cref{eq:gcp_grad_A}, the corresponding gradients can be shown to be
\begin{multline}
\label{eq:spatial_gcp_grad}
  \FPD{F}{\Ak} = \Ytk \Zk{k} \diag(\st) + \lambda\Ak + \\
  \sum_{h \in \HWt} w_h \left(\Ak \diag(\sh)\Zk{k}^T\Zk{k}\diag(\sh) - \Akold\diag(\sh)\Zkold{k}^T\Zk{k}\diag(\sh)\right) 
\end{multline}
for $k=1,\dots,d$, where 
\begin{align}
  \Zk{k}    &= \Ak{d} \odot\dots\odot \Ak{k+1} \odot \Ak{k-1} \odot\dots\odot \Ak{1}, \\
  \Zkold{k} &= \Akold{d} \odot\dots\odot \Akold{k+1} \odot \Akold{k-1} \odot\dots\odot \Akold{1}.
\end{align}

\subsection{Sampling for streaming stochastic GCP approximations}
\label{sec:streaming_gcp_sampling}
Since each $\Yt$ is in general dense, we must compute sampled approximations (denoted by $\Yt*$) for each SGD iteration.
In principle, any sampling method can be used, but in this work we employ the stratified sampling approach of~\cite{KoHo20} where for each time step $t$, the set of sampled coordinates $\I_t$ for computing $\Yt*$, is partitioned into two disjoint sets consisting of indices corresponding to nonzeros and zeros in $\Xt$.
As in~\cite{KoHo20}, we assume these sets are formed by sampling uniformly, {\em with replacement}, $p$ and $q$ times from the sets of nonzeros and zeros, respectively (zeros are sampled by searching the tensor after each candidate is computed to verify the candidate is not a nonzero, and this continues until $q$ samples have been generated).
%$\{1,\dots,\eta_t\}$ and $\{1,\dots,\omega\}$, where 
For each $i\in\I$, let $\tilde{p}_{it}$ be the number of times $i$ is selected as a nonzero and $\tilde{q}_{it}$ the number of times it is selected as a zero.
Then the entries of the sampled gradient tensor $\Yt*$ are given by
\begin{equation}
  \yit* = \left(\tilde{p}_{it}\frac{\eta_t}{p} + \tilde{q}_{it}\frac{\omega-\eta_t}{q}\right)\frac{\partial f}{\partial m}(\xit,\mit)
\end{equation}
where $\eta_t = \mbox{nnz}(\Xt)$ is the number of nonzeros in $\Xt$ and $\omega = \prod_{k=1}^d I_k$ is the total number of elements of $\Xt$ (which is independent of $t$).
Since $\E[\tilde{p}_{it}] = p/\eta_t$ and $\E[\tilde{q}_{it}] = q/(\omega-\eta_t)$, it is easy to see that $\E[\yit*] = \yit$.  The history and regularization terms in \cref{eq:spatial_gcp_grad} could also be sampled, but since their true values can be computed efficiently, there is no reason to do so and their true gradient contributions are included in each stochastic gradient.

\begin{table}
  \centering%\footnotesize 
  \caption{OnlineGCP input arguments and algorithmic parameters.}
  \label{tab:online_gcp_args}
  \begin{tabular}{|l|l|}
    \hline
    \bf Argument & \bf Description \\ \hline
     $\{\Xt\}_{t=1,2,\dots}$ & Streamed tensor slices \\
     $\{\Ak{k}\}$ & Initial guess for $k=1,\dots,d$ factor matrices \\
     $\HW$ & Initial history window (may be empty) \\
     $H$ & Maximum size of history window \\
     $w$, $\theta$ & Exponential weighting parameters for history window \\
     $\lambda$, $\mu$ & Factor matrix/temporal weight regularization parameter \\
     $f$ & GCP loss function \\
     $p'$, $q'$ & Number of stratified nonzero/zero samples for estimating $F$ \\
     $p$, $q$ & Number of stratified nonzero/zero samples for estimating $\Y$ \\
     \textsc{AdamWght} & Instance of the \textsc{Adam} class for the temporal weights solver \\
     \textsc{AdamFac} & Instance of the \textsc{Adam} class for factor matrix solver \\
     $\tol_w$, $\tol_f$ & Tolerances for weights/factor matrix solvers \\
     $\kappa_w$, $\kappa_f$ & Maximum number of epochs for weights/factor matrix solvers \\
     $\tau_w$, $\tau_f$ & Number of iterations per epoch for weights/factor matrix solvers \\
    \hline
  \end{tabular}
\end{table}

Similar sampling calculations are required for efficiently approximating the objective function $F$ in~\cref{eq:gcp_stream_objective}, however as in~\cite{KoHo20} there are a few changes.
First, we use a much larger number of samples when approximating $F$ to ensure accuracy.
Second, we use the same set of samples across all epochs within the temporal and factor matrix solvers for consistent estimations of convergence (but compute a different set of samples for the temporal and factor matrix solvers, and also for each slice $\Xt$).
As before, we sample uniformly with replacement $p'$ and $q'$ times from the sets of indices corresponding to nonzeros and zeros, respectively.
As in the gradient, the true value of the history and regularization objective terms can be efficiently computed, so the sampled approximation $\tilde{F}$ to $F$ is given by
\begin{multline}
  \tilde{F}(\Xt,\Ak{1},\dots,\Ak{d},\st) = \\
  \sum_{i \in \I_t} \left(\tilde{p}'_{it}\frac{\eta_t}{p'}+\tilde{q}'_{it}\frac{\omega-\eta_t}{q'}\right)f(\xit,\mit)
  + \sum_{h\in\HWt} w_h\|\Mhold-\Mh\|_F^2  \\
  + \frac{\lambda}{2} \sum_{k=1}^d \| \Ak \|_F^2 + \frac{\mu}{2} \| \st \|_2^2,
\end{multline}
with similar definitions of $\tilde{p}'_{it}$ and $\tilde{q}'_{it}$.  It is straightforward to see that $\E[\tilde{F}] = F$.

\begin{algorithm}[t]
\caption{Streaming GCP algorithm OnlineGCP.}
\label{alg:streaming_gcp}
\begin{algorithmic}[1]
\Function{OnlineGCP}{$\{\Xt\}_{t=1,2,\dots}$, $\{\Ak{k}\}$, $\HW$, $\{\sh\}_{h\in\HW}$}
\State \textsc{AdamFac.Init}($\{\Ak{k}\}$) \Comment{Initialize factor matrix ADAM step object} \label{line:adam_init}
\State $i \gets 0$ \Comment{Factor matrix update iteration index}
\For{$t=1,2,\dots$} \label{line:tensor_loop}
  \For{$k=1,\dots,d$}
    \State $\Ak \gets \Akold$
  \EndFor
  \State $\st \gets \textsc{WghtGcpSgd}\big(\Xt, \{\Ak\},t)$ \label{line:temp_gcp}
  \State $\{\Ak\},i \gets \textsc{FacGcpSgd}\big(\Xt, \{\Ak\}, \st, \{\Akold\}, \HW, \{\sh\}_{h\in\HW},t,i)$ \label{line:fac_gcp}
  %\For{$j=1,\dots,\tau_s$}\label{line:fac_gcp_begin}
  %  \State $i \gets i+1$
  %  \State $\begin{aligned}\{\Akt\} \gets \textsc{GcpFactorUpdate}\big(&\Xt, \{\Akt\}, \st, \{\Akold\}, \HW, \{\sh\}_{h\in\HW}, \\& w, \theta, \lambda, t, f, l, p_g, q_g, \textsc{AdamFac}, i)\end{aligned}$
  %\EndFor \label{line:fac_gcp_end}
  \If{$|\HW|<H$} \Comment{Update history window}\label{line:hist_update_begin}
    \State $\HW \gets \HW \cup \{t\}$
  \Else \Comment{Reservoir sampling to compute $\HWt{t+1}$}
    %\State $\HW \gets $ random subset of size $H$ of $\HW \cup \{t\}$
    \State $j \gets$ random element of $\{1,\dots,t\}$
    \If{$j \leq H$}
      %\State $\HW(j) = t$
      \State replace $j^{th}$ element of $\HW$ with $t$
    \EndIf
  \EndIf \label{line:hist_update_end}
  \For{$k=1,\dots,d$}
    \State $\Akold \gets \Ak$
  \EndFor
\EndFor
\State \Return $\{\Ak\}$
\EndFunction
\end{algorithmic}
\end{algorithm}

\subsection{Online GCP algorithm}
\label{sec:streaming_gcp_alg}

The high level streaming OnlineGCP algorithm is summarized in~\cref{alg:streaming_gcp}.
Descriptions of the various input arguments and algorithmic parameters are summarized in~\cref{tab:online_gcp_args}.
%\tknote{I suggest to declare certain of these parameters to simply be \emph{fixed}, maybe put a solid line in the table and move these parameters to that section. This would make the algorithm statements much less complicated. Candidates for things that are fixed include: $w$, $\theta$, $\lambda$, $\mu$, both tol's, both $\kappa$'s, both $\tau$'s, and anything else that is not going to be varied in the experiments section. Also, the ``instances'' are a bit confusing. I guess that I guess what you're trying to do by saying each instance has its own memory.}
In additional to the streamed tensor slices $\{\Xt\}$, the algorithm takes on input an initial guess for the factor matrices $\{\Ak\}$ for $k=1,\dots,d$.
In our work and the results shown in~\cref{sec:numerical}, these are computed via a relevant CP factorization (using, e.g., static CP-ALS, CP-APR, or GCP methods) on an initial set of tensor slices, which we call a warm-start.
We also specify an initial history window $\HW$ to include the warm start, but this is not required.
Here $\{\sh\}_{h\in\HW}$ is the temporal mode values corresponding to the time steps contained within the history window.
The algorithm uses the ADAM SGD update strategy~\cite{KiBa15} shown in \cref{alg:adam} in both the temporal weights and factor matrix GCP solvers.
%\tknote{This ADAM is not the same as what was used by GCP, and I think that's potentially a problem. The ``UPDATE'' should walk back the factor matrices as well as the ADAM parameters. Otherwise, it's inconsistent. FWIW, you don't have to do the walk back, but if you're going to do it, do it for both. Also, the update typically shortens the step length when the update fails. Again, if you're not going to do that, then don't mess with this updating part.}
In the latter case, the ADAM first and second moment factor matrices are tracked across all streamed tensor slices, which are initialized to zero in line~\ref{line:adam_init} of \cref{alg:streaming_gcp}.  Note that in~\cref{alg:adam}, all algebraic operations (e.g., square, division, square-root) involving the relevant factor matrices/vectors are taken component-wise.
The bulk of the algorithm begins at line~\ref{line:tensor_loop} by looping over the streamed tensor slices.
For each slice, the temporal weights $\st$ are computed using a variant of the GCP-SGD solver algorithm in line~\ref{line:temp_gcp}.
Then in line~\ref{line:fac_gcp}, the factor matrices are updated.
Finally, the history window is updated in lines~\ref{line:hist_update_begin}--\ref{line:hist_update_end}.

\algblockdefx[Class]{Class}{EndClass}%
  [1][Unknown]{\textbf{class} #1:}%
  {\textbf{end class}}
\algblockdefx[Properties]{Properties}{EndProperties}%
  {\textbf{properties:}}%
  {\textbf{end properties}}
\begin{algorithm}[t]
\caption{ADAM SGD Step Method.}
\label{alg:adam}
\begin{algorithmic}[1]
\Class[\textsc{Adam}]
  \Properties
    \State $\alpha$, $\beta_1$, $\beta_2$, $\epsilon$ \Comment{Update parameters}
    \State $\Vc{u}$, $\Vc{v}$ \Comment{First and second moments}
    \State $\Vc{u}_o$, $\Vc{v}_o$ \Comment{First and second moments from previous epoch}
    \State $\Vc{a}_o$ \Comment{Solution from previous epoch}
    \State $l$ \Comment{Lower bound}
  \EndProperties
  \Function{Init}{\Vc{a}}
    \State $\Vc{u}, \Vc{v}, \Vc{u}_o, \Vc{v}_o , \Vc{a}_o \gets \textsc{zeros}(\textsc{size}(\Vc{a}))$ \Comment{Initialize to zero with shape of $\Vc{a}$}
  \EndFunction
  \Function{Step}{$\Vc{a}$, $\Vc{g}$, $i$}
    \State $\Vc{u} \gets \beta_1\Vc{u} + (1-\beta_1)\Vc{g}$ \Comment{ADAM update for first moment $\Vc{u}$}
    \State $\Vc{v} \gets \beta_2\Vc{v} + (1-\beta_2)\Vc{g}^2$ \Comment{ADAM update for second moment $\Vc{v}$}
    \State $\tilde{\alpha} \gets \alpha\sqrt{1-\beta_2^i}/(1-\beta_1^i)$ \Comment{ADAM bias-corrected learning rate}
	\State $\Vc{a} \gets \Vc{a} - \tilde{\alpha}\Vc{u}/(\sqrt{\Vc{v}}+\epsilon)$ \Comment{ADAM update}
	\State $\Vc{a} \gets \mbox{max}(\Vc{a}, l)$ \Comment{Lower bound}
  \EndFunction
  \Function{Update}{$\Vc{a}$, passed}
  	\If{passed} \Comment{Set old moments to newest values}
		\State $\Vc{u}_o \gets \Vc{u}$ 
		\State $\Vc{v}_o \gets \Vc{v}$
		\State $\Vc{a}_o \gets \Vc{a}$
	\Else \Comment{Reset moments to saved values}
		\State $\Vc{u} \gets \Vc{u}_o$ 
		\State $\Vc{v} \gets \Vc{v}_o$
		\State $\Vc{a} \gets \Vc{a}_o$
		\State Decrease learning rate $\alpha$
	\EndIf
	\State \Return $\Vc{a}$
  \EndFunction
\EndClass
\end{algorithmic}
\end{algorithm}

\begin{algorithm}%[t]
\caption{GCP-SGD solver for temporal weights $\mathbf{s}$.}
\label{alg:temporal_gcp_sgd}
\begin{algorithmic}[1]
\Function{WghtGcpSgd}{$\Xt$, $\{\Ak\}$, $t$}
  \State $\Vc{s} \gets$ zero vector of length $R$ \Comment{initialize $\Vc{s}$}
  \State \textsc{AdamWght.init}($\Vc{s}$) \Comment{Reset ADAM stepper}
  \State $\fest \gets \textsc{EstimateObjective}(\Xt, \{\Ak\}, \Vc{s}, \{\}, \{\}, \{\}, t, 0, 0, \mu)$
  \State $i \gets 0$ \Comment{Iteration counter}
  \While{$\fest > \tol_w$ and at most $\kappa_w$ iterations}
  \State $\festold \gets \fest$
  \For{$\tau_w$ iterations}
    \State $\Y \gets \textsc{StratGrad}(\Xt, \{\Ak\}, \Vc{s})$ %\Comment{GCP sample tensor}
    \State $\Vc{g} \gets \Z'\mathbf{y} + \lambda\Vc{s}$ \Comment{GCP gradient for temporal weights $\Vc{s}$} \label{line:gcp_temp_grad}
    \State $\Vc{s} \gets$ \textsc{AdamWght.Step}($\Vc{s}$, $\Vc{g}$, $i$) \Comment{Apply ADAM step}
	\State $i \gets i+1$
  \EndFor
  \State $\fest \gets \textsc{EstimateObjective}(\Xt, \{\Ak\}, \Vc{s}, \{\}, \{\}, \{\}, t, 0, 0, \mu)$  %\Comment{estimate loss with fixed set of samples}
  \If{$\fest > \festold$} %\Comment{check for failure to decrease loss}
  	\State $\Vc{s}\gets\textsc{Adam.update}(\Vc{s},\textsc{false})$
    \State $\fest \gets \festold$  %\Comment{revert to last epoch's variables}
    \State $i \gets i - \tau_w$
  \Else
  	\State $\Vc{s}\gets\textsc{AdamWght.update}(\Vc{s},\textsc{true})$
  \EndIf
  \EndWhile
\State \Return $\Vc{s}$
\EndFunction
\end{algorithmic}
\end{algorithm}

\begin{algorithm}[t]
\caption{GCP-SGD solver for factor matrices $\{\Ak\}$.}
\label{alg:gcp_factor_update}
\begin{algorithmic}[1]
\Function{FacGcpSgd}{$\Xt$, $\{\Ak\}$, $\st$, $\{\Akold\}$, $\HW$, $\{\sh\}_{h\in\HW}$, $t$, $i$}
  \State $\fest \gets \textsc{EstimateObjective}(\Xt, \{\Ak\}, \st, \{\Akold\}, \HW, \{\sh\}_{h\in\HW}, t, w, \lambda, 0)$
  \While{$\fest > \tol_f$ and at most $\kappa_f$ iterations}
  \State $\festold \gets \fest$
  \For{$\tau_f$ iterations}
    \State $\Y \gets \textsc{StratGrad}(\Xt, \{\Ak\}, \st)$ %\Comment{GCP sample tensor}
    \For{$k=1,\dots,d$} \label{line:gcp_fac_grad_beg} %\Comment{GCP gradient for factor matrix $\Ak$}
   	  \State $\Gk \gets \Yk \Zk \diag(\st) + \lambda\Ak$ 
      \For{$h\in\HW$}
        \State $\Mx{B} \gets \diag(\sh)\Zk{k}^T\Zk{k}\diag(\sh)$
        \State $\Mx{C} \gets \diag(\sh)\Zkold{k}^T\Zk{k}\diag(\sh)$
        \State $\Gk \gets \Gk + w \theta^{t-h} \left(\Ak\Mx{B}  - \Akold\Mx{C}\right) $ 
      \EndFor
    \EndFor \label{line:gcp_fac_grad_end}
    \State $\{\Ak\} \gets$ \textsc{AdamFac.Step}($\{\Ak\}$, $\{\Gk\}$, $i$) %\Comment{Apply ADAM step}
	\State $i \gets i+1$
  \EndFor
  \State $\fest \gets \textsc{EstimateObjective}(\Xt, \{\Ak\}, \st, \{\Akold\}, \HW, \{\sh\}_{h\in\HW}, t, w, \lambda, 0)$
  \If{$\fest > \festold$} %\Comment{check for failure to decrease loss}
  	\State $\{\Ak\}\gets\textsc{AdamFac.update}(\{\Ak\},\textsc{false})$
    \State $\fest \gets \festold$  %\Comment{revert to last epoch's variables}
    \State $i \gets i - \tau_f$
  \Else
  	\State $\{\Ak\}\gets\textsc{Adam.update}(\{\Ak\},\textsc{true})$
  \EndIf
  \EndWhile
\State \Return $\{\Ak\}$, $i$
\EndFunction
\end{algorithmic}
\end{algorithm}

\begin{algorithm}%[t]
\caption{Objective estimate using stratified sampling with $p'$ nonzeros of $\X$ and $q'$ zeros.}
\label{alg:f_est}
\begin{algorithmic}[1]
\Function{EstimateObjective}{$\X$, $\{\Ak\}$, $\st$, $\{\Akold\}$, $\HW$, $\{\sh\}_{h\in\HW}$, $t$, $w$, $\lambda$, $\mu$}
\State $\eta \gets \mbox{nnz}(\X)$ \Comment{Number of nonzeros in $\X$}
\State $\zeta \gets \prod_k I_k - \eta$ \Comment{Number of zeros in $\X$}
\State $\fest \gets 0$
\For{$c=1,\dots,p'$} \Comment{Sample $p$ nonzeros of $\X$}
  \State $i \gets$ random element of $\{1,\dots,\eta\}$
  \State $(i_1,\dots,i_d) \gets$ indices of nonzero $i$ of $\X$
  \State $m \gets \sum_{j=1}^R s_{tj} \prod_{k=1}^d a^{(k)}_{i_k j}$
  \State $\fest \gets \fest + (\eta/p')f(x_{i_1,\dots,i_p},m)$
\EndFor
\State $c\gets 1$
\While{$c\leq q'$} \Comment{Sample $q$ zeros of $\X$}
  \For{$k=1,\dots,d$}
    \State $i_k \gets$ random element of $\{1,\dots,I_k\}$
  \EndFor
  \If{$x_{i_1,\dots,i_d} \neq 0$}
    \State \textbf{continue}
  \EndIf
  \State $m \gets \sum_{j=1}^R s_{tj} \prod_{k=1}^d a^{(k)}_{i_k j}$
  \State $\fest \gets \fest + (\zeta/q')f(0,m)$
  \State $c\gets c+1$
\EndWhile
\State \Return $\fest + \frac{w}{2} \sum_{h\in\HWt} \theta^{t-h} \|\Mhold-\Mh\|_F^2 + \frac{\lambda}{2} \sum_{k=1}^d \| \Ak \|_F^2 + \frac{\mu}{2} \|\st\|_2^2$
\EndFunction
\end{algorithmic}
\end{algorithm}

\begin{algorithm}%[t]
\caption{Stratified sampling for $p$ nonzeros of $\X$ and $q$ zeros.}
\label{alg:strat_grad_sampling}
\begin{algorithmic}[1]
\Function{StratGrad}{$\X$, $\{\Ak\}$, $\st$}
\State $\eta \gets \mbox{nnz}(\X)$ \Comment{Number of nonzeros in $\X$}
\State $\zeta \gets \prod_k I_k - \eta$ \Comment{Number of zeros in $\X$}
\State $\tilde{\Y} \gets 0$
\For{$c=1,\dots,p$} \Comment{Sample $p$ nonzeros of $\X$}
  \State $i \gets$ random element of $\{1,\dots,\eta\}$
  \State $(i_1,\dots,i_d) \gets$ indices of nonzero $i$ of $\X$
  \State $m \gets \sum_{j=1}^R s_{tj} \prod_{k=1}^d a^{(k)}_{i_k j}$
  \State $\tilde{y}_{i_1,\dots,i_d} \gets \tilde{y}_{i_1,\dots,i_d} + (\eta/p)(\partial f/\partial m)(x_{i_1,\dots,i_p},m)$
\EndFor
\State $c\gets 1$
\While{$c\leq q$} \Comment{Sample $q$ zeros of $\X$}
  \For{$k=1,\dots,d$}
    \State $i_k \gets$ random element of $\{1,\dots,I_k\}$
  \EndFor
  \If{$x_{i_1,\dots,i_d} \neq 0$}
    \State \textbf{continue}
  \EndIf
  \State $m \gets \sum_{j=1}^R s_{tj} \prod_{k=1}^d a^{(k)}_{i_k j}$
  \State $\tilde{y}_{i_1,\dots,i_d} \gets \tilde{y}_{i_1,\dots,i_d} + (\zeta/q)(\partial f/\partial m)(0,m)$
  \State $c\gets c+1$
\EndWhile
\State \Return $\tilde{\Y}$
\EndFunction
\end{algorithmic}
\end{algorithm}

The temporal weights and factor matrix GCP-SGD solver algorithms are summarized in~\cref{alg:temporal_gcp_sgd} and~\cref{alg:gcp_factor_update}, which are in general similar to the corresponding algorithm in~\cite{KoHo20} except modified for solving for just the temporal weights or factor matrices, respectively.
As described above, they uses the stratified sampling approach of~\cite{KoHo20} for estimating the objective function $F$, which is summarized in~\cref{alg:f_est}, and for estimating the gradient summarized in~\cref{alg:strat_grad_sampling}.
Note that since the history window does not affect the gradient for the current temporal weights, it is not included in the objective function and gradient calculations in the temporal solver.

%%% Local Variables:
%%% mode: latex
%%% TeX-master: "online_gcp_paper"
%%% End:

\section{Software Implementation}
\label{sec:software}
% !TEX root = online_gcp_paper.tex

Our software implementation of \cref{alg:streaming_gcp,alg:adam,alg:temporal_gcp_sgd,alg:gcp_factor_update,alg:f_est,alg:strat_grad_sampling} is written in a hybrid of Matlab and C++ code leveraging the Matlab Tensor Toolbox~\cite{BaKo07,TensorToolbox} and the C++ GenTen package for performance portable tensor decompositions~\cite{GentenSISC,GenTen}.  Most of the code is written in object oriented fashion for Matlab leveraging the Tensor Toolbox for representing sparse tensors and K-tensors.  In particular, \cref{alg:streaming_gcp} is implemented in Matlab as a function that takes as an argument all of the tensor slices to be analyzed, looping over them as in the for-loop starting at line~\ref{line:tensor_loop}.  The contents of the for-loop are then implemented through a separate Matlab class that computes the new temporal weights and updates the factor matrices for each tensor slice.  The purpose of this design was to allow the contents of the loop to be executed in a true streaming context where not all tensor slices are available beforehand.  The temporal weights SGD solver (\cref{alg:temporal_gcp_sgd}) and factor matrix solver (\cref{alg:gcp_factor_update}) are encapsulated in a unified class hierarchy implementing the general GCP-SGD solver strategy for a single sparse tensor.  Arguments are supplied controlling which modes/weights are to be updated.  In the temporal case, only the temporal weights are updated while the factor matrices are held fixed.  Similarly, the factor matrix updates keep the temporal weights fixed.  Thus the code implements a static GCP factorization procedure as a special case.  

This class hierarchy abstracts the sampling procedures for approximating $F$ and $\Y$ (e.g., stratified) and the SGD step procedures (e.g., ADAM) in~\cref{alg:adam,alg:temporal_gcp_sgd,alg:gcp_factor_update} allowing new sampling and step procedures to be easily plugged in.  Furthermore, the sampling classes also implement the needed gradient computations for the GCP-SGD procedure (i.e., line~\ref{line:gcp_temp_grad} of~\cref{alg:temporal_gcp_sgd} and lines~\ref{line:gcp_fac_grad_beg}--\ref{line:gcp_fac_grad_end} of~\cref{alg:gcp_factor_update}), allowing the code to take advantage of special structure.  For example, a ``dense'' sampling class is provided for the Gaussian case with $f(x,m) = (x-m)^2$ that uses for the factor matrix gradient (not including regularization and history terms for brevity)
\begin{equation}
  \frac{1}{2}\FPD{F}{\Ak} = \Xk \Zk{k} \diag(\st)-\Ak\diag(\st)\Zk{k}'\Zk{k}\diag(\st).
\end{equation}
The code also supports replacing the GCP-SGD solution procedure for the temporal weights with a single least-squares solve in the Gaussian case, and thus includes the original OnlineSGD algorithm as a special case by choosing the Gaussian loss function, least-squares temporal solve, dense gradient evaluation, no history window, and a single factor matrix update iteration per streamed tensor slice.

Finally, the sampling/gradient and step abstractions allow for plugging in high-performance, multi-threaded implementations of these computations provided by the GenTen software package.  This software package implements sampling procedures, MTTKRP, and SGD step procedures using the Kokkos C++ performance portability API~\cite{Kokkos:2012:SciProg,Kokkos:2014:JPDC} allowing a single C++ implementation of each kernel to be executed with high performance on a variety of contemporary architectures, including multicore CPUs and manycore GPUs.  These kernels are exposed to Matlab through the MEX API and a variety of bundled Matlab classes and functions allowing GenTen kernels to be integrated within the Tensor Toolbox.  This enables the high-level streaming algorithm to be rapidly developed and modified within Matlab but enable high-performance by executing all performance-critical kernels with compiled C++, multithreaded code.  For example, this allows the streaming code to be executed in Matlab, but have the kernels run on attached GPUs, resulting in substantial speedups.
%%% Local Variables:
%%% mode: latex
%%% TeX-master: "online_gcp_paper"
%%% End:

\section{Numerical Experiments}
\label{sec:numerical}
% !TEX root = online_gcp_paper.tex

We now present several numerical experiments that compare the accuracy of the OnlineGCP method with several static and streaming alternatives for Gaussian, Poisson, and Bernoulli loss functions.  The static methods are applied to  the entire $d+1$-way tensor formed by stacking the streamed slices across the temporal mode, whereas the streaming methods update their decomposition one slice at a time.  For the streaming methods, we generate an initial CP model by applying an appropriate static decomposition method (i.e., CP-ALS, CP-APR, or GCP) to a small portion of the streaming data, which we call a warm-start.
Throughout these experiments we will measure the effectiveness of the OnlineGCP approach by comparing $\it{local}$ and $\it{global}$ reconstruction losses.
In the context of streaming we define the local reconstruction loss as the total loss for a given decomposition, divided by the norm of the data, for every observed time slice:
\begin{equation}
\label{eq:local_loss}
  F_{local}(\Xt,\Mt) = \frac{1}{\|\X_t\|_F^2}\sum_{i\in\I} f(\xit,\mit)
\end{equation}
where $\Mt=\KTt$ and $\Ak{1},\dots,\Ak{d}$ indicate the values of factor matrices computed at that point in time.
This evaluates how well our current model fits the most recently observed data.
For OnlineGCP, we compute a sampled approximation to $F_{local}$ using the sampling procedure in \cref{alg:f_est} but not including history, factor matrix, or temporal weight regularization terms.

For the global reconstruction loss we back test the model at our final time step against all previously observed data.  The functional form is the same as in \cref{eq:local_loss} however we use the factor matrices $\Ak{1},\dots,\Ak{d}$ from the final time step.  This evaluates how well the final model approximates all observed data.  For OnlineGCP, we compute the true value of the global loss instead of a sampled approximation.  Since the static methods use all available data, the local and global reconstruction losses are identical.

As indicated in \cref{tab:online_gcp_args}, the OnlineSGD method has numerous hyperparameters.  The values used in each experiment are summarized in \cref{tab:hyperparameters}.  These values were chosen empirically to produce good results through hand-tuning, but are not necessarily optimal.
\begin{table}[ht]
\centering
\caption{OnlineSGD hyperparameters for the numerical experiments.  Here $R$ is the rank of the CP model, $\alpha_w$, $\alpha_f$ are the ADAM learning rates of the temporal weight and factor matrix solvers, respectively, $\kappa_w$, $\kappa_f$ are the number of epochs for each solver, $w$ is the multiplicative weight of the history term, $H$ is the size of the history window, and $H_{init}$ is the size of the warm-start.  All experiments used $\theta=1$ (no exponential down-weighting of historical slices), $\lambda=\mu=0$ (no rank regularization penalty), $\tau_w=\tau_f=100$ iterations per epoch, and the default ADAM update parameters ($\beta_1=0.9$, $\beta_2=0.999$, and $\epsilon=10^{-8}$).}
\label{tab:hyperparameters}
\begin{tabular}{crrrcrrrc}
  \toprule
  {\bf Experiment}       & $R$ & $\alpha_w$ & $\kappa_w$ & $\alpha_f$      & $\kappa_f$ & $w$ & $H$ & $H_{init}$ \\
  \midrule
  { Synthetic Gaussian}  & 20  & 10.0       & 20         & $1\cdot10^{-4}$ & 5          & 1   & 50  & 10 \\
  { Synthetic Poisson}   & 20  & 1.0        & 20         & $1\cdot10^{-4}$ & 10         & 10  & 50  & 10 \\
  { Taxicab Poisson}     & 50  & 10.0       & 1          & $1\cdot10^{-3}$ & 1          & 1   & 30  & 20 \\
  %{ ArXiv Poisson}       & 50  & 10.0       & 5          & $5\cdot10^{-4}$ & 10         & 1   & 50  & 10 \\
  { Chicago Binary}      & 50  & 0.1        & 5          & $1\cdot10^{-3}$ & 5          & 10  & 500 & 20 \\
  \bottomrule
\end{tabular}
\end{table}
\begin{figure*}[h!]
	\centering
	\includegraphics[width=\textwidth]{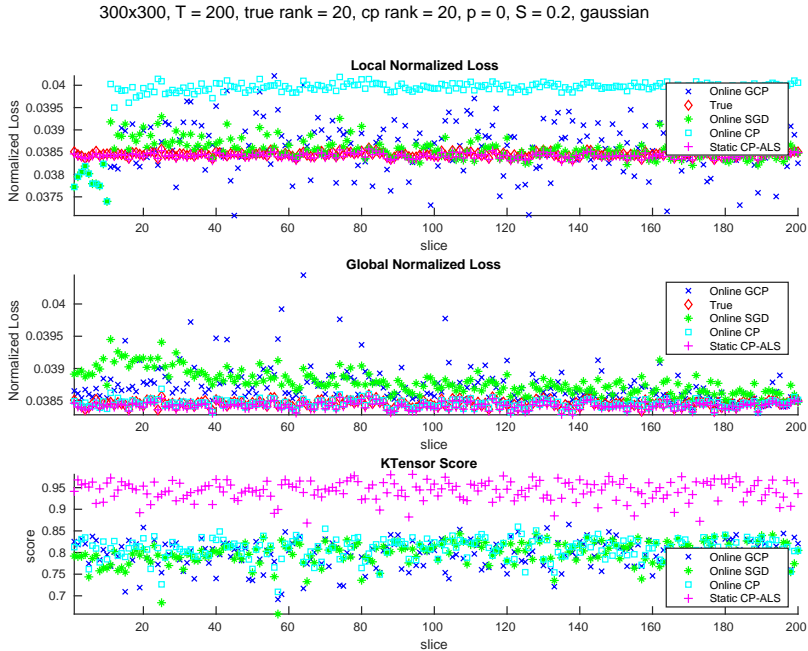}
	\caption{{\it Results of the synthetic Gaussian experiment showing comparable performance between OnlineGCP, OnlineSGD, OnlineCP, and static CP-ALS.}}
	\label{fig:synthetic_gaussian}
\end{figure*}

\subsection{Synthetic Data Experiments}
\label{sec:synthetic}

We first describe experiments with two synthetic data sets derived from randomly generated K-tensors which the computed decomposition methods should recover.  For these experiments we also measure the congruence~\cite{TOMASI20061700} between the K-tensor computed by each method and K-tensor used to generate the data (called the K-tensor score).  A perfect recovery corresponds to a score of 1.0.

\begin{figure*}
	\centering
	\includegraphics[width=\textwidth]{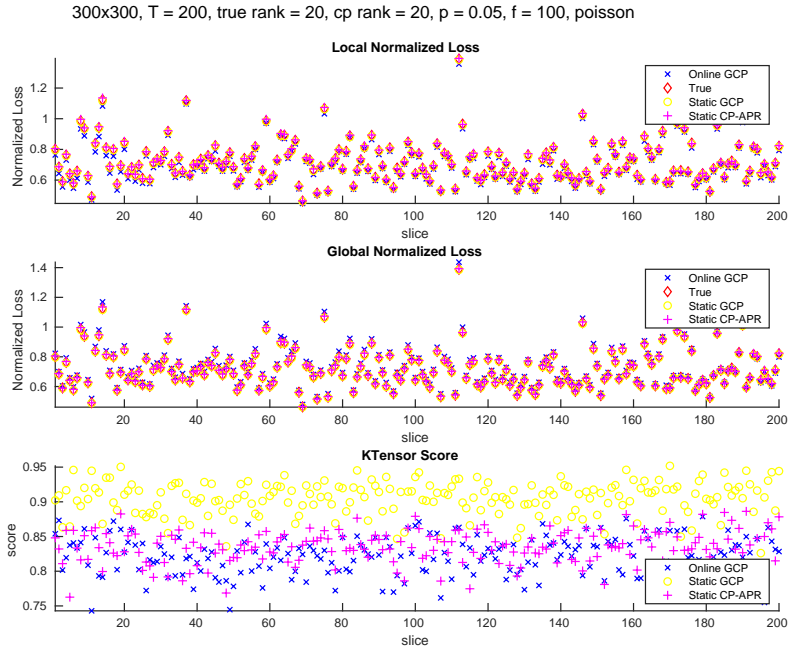}
	\caption{{\it Results of the synthetic Poisson experiment showing comparable performance between OnlineGCP, static CP-APR, and static GCP (using 50,000 zero/nonzero samples for each objective and gradient evaluation).}}
	\label{fig:synthetic_poisson}
\end{figure*}

\subsubsection*{Gaussian}
To construct a Gaussian distributed synthetic data set, we first constructed a random 3-way, rank-20 K-tensor with factor matrices of size $300\times20$, $300\times20$, and $200\times20$, respectively.  Each factor matrix entry was drawn uniformly at random from $(0,1)$.  This K-tensor provides the ground truth for the model.  A dense tensor was then generated by multiplying out the K-tensor and perturbing each entry by draws from a zero-mean Gaussian distribution with a standard deviation of 0.2.  This tensor is then streamed slice-by-slice, where each slice is a dense $300\times 300$ matrix.  The warm-start was generated by applying CP-ALS to the first 10 slices.  
We then compared our method to static CP-ALS applied to the full $300\times300\times200$ tensor, OnlineCP, and OnlineSGD using 10,000 nonzero samples for each objective/gradient evaluation, no zero samples (since the tensor is dense), and the remaining hyperparameters as indicated in \cref{tab:hyperparameters}.
In \cref{fig:synthetic_gaussian} we demonstrate comparable results for the local reconstruction loss, global reconstruction loss, and K-tensor score with respect to the ground truth and the considered methods.

\subsubsection*{Poisson}
To generate a Poisson distributed synthetic data set, we used the procedure described in \cite{ChKo12} to generate a sparse 3-way tensor of size $300\times300\times200$, $R=20$ factors, and roughly 3.2\% nonzero sparsity.  We then stream this tensor slice-by-slice in OnlineGCP as before, but this time comparing to static CP-APR and GCP with Poisson loss computed from the full tensor dataset.  OnlineGCP used a warm-start constructed from applying CP-APR to the first 10 slices, all tensor nonzeros in the objective/gradient evaluations, and 50,000 and 10,000 zero samples for the objective and gradient, respectively.
In Figure~\ref{fig:synthetic_poisson} we see fairly comparable results in losses among all of the methods.
For the comparison to the ground truth decomposition CP-APR performs slightly better than the two GCP based methods.

\subsection{Realistic Data Experiments}
\label{sec:real}
We now present several experiments using real data tensors with non-Gaussian loss functions.

\begin{figure*}
	\centering
	\includegraphics[width=\textwidth]{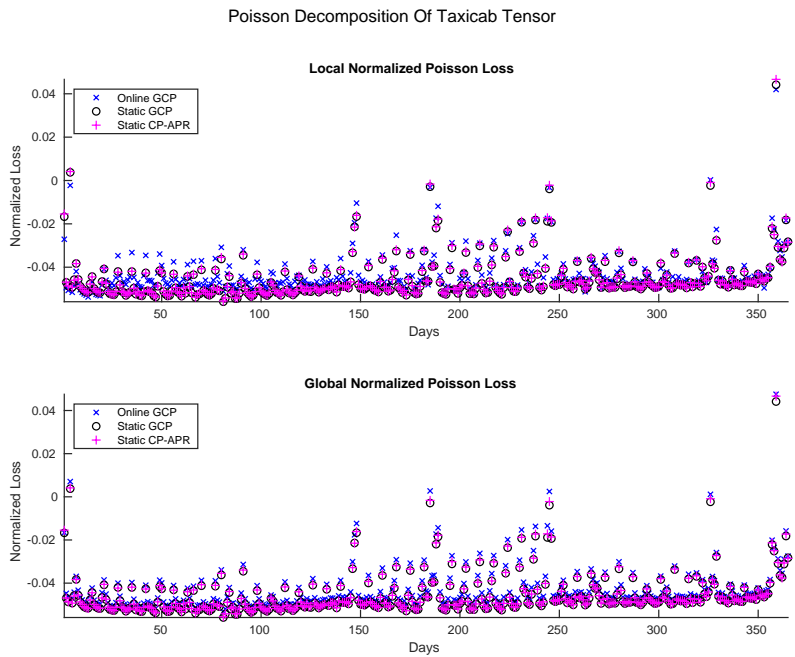}
	\caption{{\it Results of the NYC Taxicab experiment with Poisson loss showing comparable performance between OnlineGCP, static CP-APR, and static GCP (using 50,000 zero/nonzero samples for each objective evaluation and 10,000 zero/nonzero samples for each gradient evaluation).}}
	\label{fig:taxi_poisson}
\end{figure*}

\subsubsection*{NYC Taxicab}
We use New York City Yellow Taxi data from 2018 to generate tensors corresponding to travel throughout the city \cite{TaxiData}.
These data provide fields containing pick-up and drop-off dates/times, pick-up and drop-off locations, trip distances, itemized fares, rate types, payment types, and driver-reported passenger counts.
For this experiment, we use pick-up and drop-off dates/times to generate a four way tensor of size $263\times263\times24\times365$ with approximately 3.8\% nonzero sparsity containing counts of the number of taxi rides between the given taxi zones over the given hour for the given day.  Accordingly, we generated CP decompositions using Poisson loss via OnlineGCP (streaming one slice at a time corresponding to a single day), GCP, and CP-APR.  OnlineGCP used a warm-start constructed from applying CP-APR to the first 10 days worth of data, 50,000 zero and nonzero samples for the objective, and 10,000 zero and nonzero samples for the gradient.
In Figure~\ref{fig:taxi_poisson} we again see comparable results in terms of the achieved local/global loss compared to the static methods.

\begin{figure*}
	\centering
	\includegraphics[width=\textwidth]{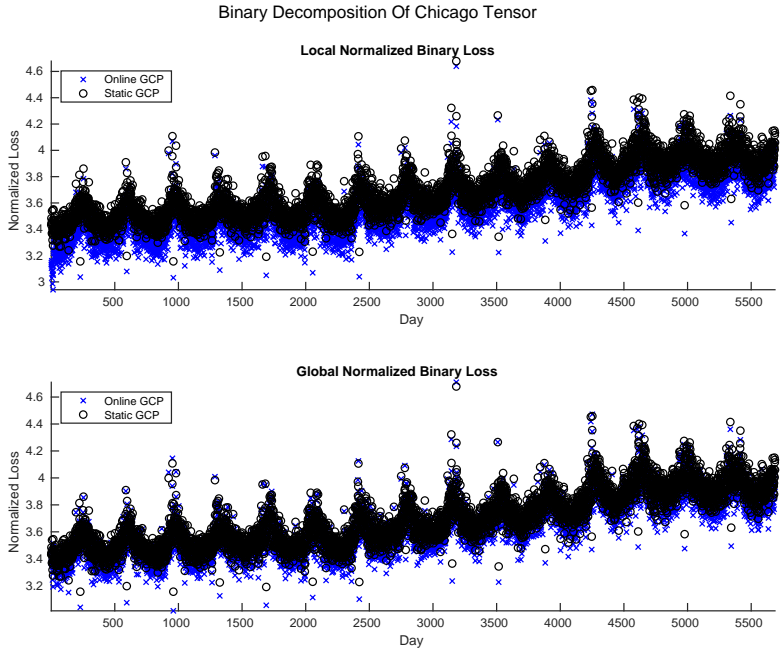}
	\caption{{\it Results of the Chicago Crime experiment with Bernoulli loss showing comparable performance between OnlineGCP and static GCP (using 50,000 zero/nonzero samples for each objective evaluation and 10,000 zero/nonzero samples for each gradient evaluation).}}
	\label{fig:chicago_binary}
\end{figure*}

\subsubsection*{Chicago Crime}
To demonstrate the approach for non-Poisson loss, we used the Chicago Crime tensor provided by FROSTT \cite{frosttdataset} converted to a binary tensor where any nonzero value was replaced by one (this is reasonable since a majority of the entries are one anyway).  
Sticking with our streaming convention across the last mode we oriented the tensor such that entry $\X(i,j,k,l)$ denoted whether on hour $i$, crime $j$ was committed in neighborhood $k$ for day $l$ from our first date.  We used a starting data of $l=500$ because there is significantly less data for days prior to this date.  The resulting tensor is of size $24\times77\times32\times5687$ with roughly 1.6\% sparsity.  A warm-start for the first 20 days was generated via GCP with Bernoulli loss using 50,000 zero/nonzero samples for the objective function and 10,0000  zero/nonzero samples for the gradient.  Given the relatively small tensor slices each time step, OnlineGCP used all nonzeros along with 10,000 and 1,000 zero samples for each objective function and gradient evaluation, respectively.  In Figure~\ref{fig:chicago_binary} we again see comparable results in terms of the achieved local/global loss compared to the static GCP method.  Note however that OnlineGCP required a much larger history window ($H=500$) than the other experiments to maintain consistent global loss over the entire streaming experiment.

\section{Conclusions}
\label{sec:conclusions}
% !TEX root = online_gcp_paper.tex
In this work, we developed a method called OnlineGCP for efficiently computing GCP decompositions of streaming tensor data.  The method extends prior work in the literature on streaming CP decompositions to the GCP case allowing for arbitrary objective/loss functions defining the CP optimization problem.  Similar to other streaming CP methods, the approach incrementally updates the temporal weights and CP model factors as each new tensor slice is observed without revisiting prior data.  It includes a tunable history term to balance reconstruction of new and old tensor data, and employs stochastic gradient descent solvers enabling scalability to large, sparse tensors.  The effectiveness of the approach was demonstrated on several synthetic and real datasets incorporating Guassian, Poisson, and Bernoulli loss functions, where comparable losses were observed compared to other streaming and static methods appropriate for the chosen form of loss.

While the approach was shown to be effective and scalable, it relies on expert choice of numerous hyperparameters that can dramatically affect accuracy and computational cost.  Unfortunately, our experience has shown these parameters must be empirically chosen on a case-by-case basis.  The sensitivity of the method to these hyperparameters primarily derives from its use of stochastic gradient descent as an optimization strategy, so future work will involve investigation of alternative solution strategies that rely on fewer hyperparameters and are more robust to their values.

\section*{Acknowledgments}
\label{sec:ack}
We would like to thank Kyle Gilman for numerous insightful comments regarding the formulation of the streaming GCP problem.

%\appendix
%\section{Arxiv History}
%\label{sec:arxiv-drift}
%\input{sec-arxiv-drift}

%\section{\LaTeX\ Macros}
%\label{sec:macros}
%\input{sec-macros}

\bibliographystyle{siamplain}
\bibliography{references,bibfiles/allrefs,bibfiles/tgk}

\end{document}